\documentclass[11pt]{article}
\usepackage{amsmath}
\usepackage{amssymb}
\usepackage{amsthm}
\usepackage{mathrsfs}
\begin{document}
\title{The boundedness of some operators with rough kernel on the weighted Morrey spaces}
\author{Hua Wang\,\footnote{E-mail address: wanghua@pku.edu.cn.}\\[3mm]
\footnotesize{School of Mathematical Sciences, Peking University, Beijing 100871, China}}
\date{}
\maketitle
\begin{abstract}
Let $\Omega\in L^q(S^{n-1})$ with $1<q\le\infty$ be homogeneous of degree zero and has mean value zero on $S^{n-1}$. In this paper, we will study the boundedness of homogeneous singular integrals and Marcinkiewicz integrals with rough kernel on the weighted Morrey spaces $L^{p,\kappa}(w)$ for $q'\le p<\infty$(or $q'<p<\infty$) and $0<\kappa<1$. We will also prove that the commutator operators formed by a $BMO(\mathbb R^n)$ function $b(x)$ and these rough operators are bounded on the weighted Morrey spaces $L^{p,\kappa}(w)$ for $q'<p<\infty$ and $0<\kappa<1$.
\end{abstract}
\textbf{\large{1. Introduction}}
\par
Suppose that $S^{n-1}$ is the unit sphere in $\mathbb R^n$($n\ge2$) equipped with the normalized Lebesgue measure $d\sigma$. Let $\Omega\in L^q(S^{n-1})$ with $1<q\le\infty$ be homogeneous of degree zero and satisfy the cancellation condition
$$\int_{S^{n-1}}\Omega(x')d\sigma(x')=0,$$
where $x'=x/{|x|}$ for any $x\neq0$. The homogeneous singular integral operator $T_\Omega$ is defined by
$$T_\Omega f(x)=\lim_{\varepsilon\to0}\int_{|y|>\varepsilon}\frac{\Omega(y')}{|y|^n}f(x-y)\,dy$$
and a related maximal operator $M_\Omega$ is defined by
$$M_\Omega f(x)=\sup_{r>0}\frac{1}{r^n}\int_{|y|<r}|\Omega(y')f(x-y)|\,dy.$$
Let $b$ be a locally integrable function on $\mathbb R^n$, the commutator of $b$ and $T_\Omega$ is defined as follows
$$[b,T_\Omega]f(x)=b(x)T_\Omega f(x)-T_\Omega(bf)(x).$$
The Marcinkiewicz integral of higher dimension $\mu_\Omega$ is defined by
$$\mu_\Omega(f)(x)=\left(\int_0^\infty\big|F_{\Omega,t}(x)\big|^2\frac{dt}{t^3}\right)^{1/2},$$
where
$$F_{\Omega,t}(x)=\int_{|x-y|\le t}\frac{\Omega(x-y)}{|x-y|^{n-1}}f(y)\,dy.$$
It is well known that the Littlewood-Paley $g$-function is a very important tool in harmonic analysis and the Marcinkiewicz integral is essentially a Littlewood-Paley $g$-function. In this paper, we will also consider the commutator $[b,\mu_\Omega]$ which is given by the following expression
$$[b,\mu_\Omega]f(x)=\left(\int_0^\infty\big|F_{\Omega,t}^b(x)\big|^2\frac{dt}{t^3}\right)^{1/2},$$
where
$$F_{\Omega,t}^b(x)=\int_{|x-y|\le t}\frac{\Omega(x-y)}{|x-y|^{n-1}}[b(x)-b(y)]f(y)\,dy.$$
\par
The classical Morrey spaces $\mathcal L^{p,\lambda}$ were first introduced by Morrey in [10] to study the local behavior of solutions to second order elliptic partial differential equations. Recently, Komori and Shirai [9] considered the weighted version of Morrey spaces $L^{p,\kappa}(w)$ and studied the boundedness of some classical operators such as the Hardy-Littlewood maximal operator, the Calder\'on-Zygmund operator on these spaces.
\par
The main purpose of this paper is to discuss the weighted boundedness of the above operators $M_\Omega$, $T_\Omega$ and $\mu_\Omega$ with rough kernels on the weighted Morrey spaces $L^{p,\kappa}(w)$ for $q'\le p<\infty$ and $0<\kappa<1$, where we set the notation $q'=q/{(q-1)}$ when $1<q<\infty$ and $q'=1$ when $q=\infty$. We shall also show that the commutators $[b,T_\Omega]$ and $[b,\mu_\Omega]$ are bounded operators on the weighted Morrey spaces $L^{p,\kappa}(w)$ for $q'<p<\infty$ and $0<\kappa<1$, where the symbol $b$ belongs to $BMO(\mathbb R^n)$. Our main results are stated as follows.
\newtheorem{thm1}{Theorem}
\begin{thm1}
Assume that $\Omega\in L^q(S^{n-1})$ with $1<q<\infty$. Then for every $q'\le p<\infty$, $w\in A_{p/{q'}}$ and $0<\kappa<1$, there exists a constant $C>0$ independent of $f$ such that
$$\|M_\Omega(f)\|_{L^{p,\kappa}(w)}\le C\|f\|_{L^{p,\kappa}(w)}.$$
\end{thm1}
\begin{thm1}
Assume that $\Omega\in L^q(S^{n-1})$ with $1<q<\infty$. Then for every $q'\le p<\infty$, $w\in A_{p/{q'}}$ and $0<\kappa<1$, there exists a constant $C>0$ independent of $f$ such that
$$\|T_\Omega(f)\|_{L^{p,\kappa}(w)}\le C\|f\|_{L^{p,\kappa}(w)}.$$
\end{thm1}
\begin{thm1}
Assume that $\Omega\in L^q(S^{n-1})$ with $1<q<\infty$ and $b\in BMO(\mathbb R^n)$. Then for every $q'<p<\infty$, $w\in A_{p/{q'}}$ and $0<\kappa<1$, there exists a constant $C>0$ independent of $f$ such that
$$\big\|[b,T_\Omega](f)\big\|_{L^{p,\kappa}(w)}\le C\|f\|_{L^{p,\kappa}(w)}.$$
\end{thm1}
\begin{thm1}
Assume that $\Omega\in L^q(S^{n-1})$ with $1<q\le\infty$. Then for every $q'<p<\infty$, $w\in A_{p/{q'}}$ and $0<\kappa<1$, there exists a constant $C>0$ independent of $f$ such that
$$\|\mu_\Omega(f)\|_{L^{p,\kappa}(w)}\le C\|f\|_{L^{p,\kappa}(w)}.$$
\end{thm1}
\begin{thm1}
Assume that $\Omega\in L^q(S^{n-1})$ with $1<q\le\infty$ and $b\in BMO(\mathbb R^n)$. Then for every $q'<p<\infty$, $w\in A_{p/{q'}}$ and $0<\kappa<1$, there exists a constant $C>0$ independent of $f$ such that
$$\big\|[b,\mu_\Omega](f)\big\|_{L^{p,\kappa}(w)}\le C\|f\|_{L^{p,\kappa}(w)}.$$
\end{thm1}
\noindent\textbf{\large{2. Notations and definitions}}
\par
First let us recall some standard definitions and notations. The classical $A_p$ weight
theory was first introduced by Muckenhoupt in the study of weighted
$L^p$ boundedness of Hardy-Littlewood maximal functions in [11].
A weight $w$ is a locally integrable function on $\mathbb R^n$ which takes values in $(0,\infty)$ almost everywhere, $B=B(x_0,r)$ denotes the ball with the center $x_0$ and radius $r$.
We say that $w\in A_p$,\,$1<p<\infty$, if
$$\left(\frac1{|B|}\int_B w(x)\,dx\right)\left(\frac1{|B|}\int_B w(x)^{-\frac{1}{p-1}}\,dx\right)^{p-1}\le C \quad\mbox{for every ball}\; B\subseteq \mathbb
R^n,$$ where $C$ is a positive constant which is independent of $B$.\\
For the case $p=1$, $w\in A_1$, if
$$\frac1{|B|}\int_B w(x)\,dx\le C\,\underset{x\in B}{\mbox{ess\,inf}}\,w(x)\quad\mbox{for every ball}\;B\subseteq\mathbb R^n.$$
\par
A weight function $w$ is said to belong to the reverse H\"{o}lder class $RH_r$ if there exist two constants $r>1$ and $C>0$ such that the following reverse H\"{o}lder inequality holds
$$\left(\frac{1}{|B|}\int_B w(x)^r\,dx\right)^{1/r}\le C\left(\frac{1}{|B|}\int_B w(x)\,dx\right)\quad\mbox{for every ball}\; B\subseteq \mathbb R^n.$$
\par
It is well known that if $w\in A_p$ with $1<p<\infty$, then $w\in A_r$ for all $r>p$, and $w\in A_q$ for some $1<q<p$. If $w\in A_p$ with $1\le p<\infty$, then there exists $r>1$ such that $w\in RH_r$.
\par
We give the following results that we will use frequently in the sequel.
\newtheorem*{lemmaA}{Lemma A}
\begin{lemmaA}[{[6]}]
Let $w\in A_p$, $p\ge1$. Then, for any ball $B$, there exists an absolute constant $C$ such that
$$w(2B)\le C w(B).$$
In general, for any $\lambda>1$, we have
$$w(\lambda B)\le C\lambda^{np}w(B),$$
where $C$ does not depend on $B$ nor on $\lambda$.
\end{lemmaA}
\newtheorem*{lemmab}{Lemma B}
\begin{lemmab}[{[7]}]
Let $w\in RH_r$ with $r>1$. Then there exists a constant $C$ such that
$$\frac{w(E)}{w(B)}\le C\left(\frac{|E|}{|B|}\right)^{(r-1)/r}$$
for any measurable subset $E$ of a ball $B$.
\end{lemmab}
A locally integrable function $b$ is said to be in $BMO(\mathbb R^n)$ if
$$\|b\|_*=\underset{B}{\sup}\frac{1}{|B|}\int_B|b(x)-b_B|\,dx<\infty,$$
where $b_B=\frac{1}{|B|}\int_B b(y)\,dy$ and the supremum is taken over all balls $B$ in $\mathbb R^n.$
\newtheorem*{thmc}{Theorem C}
\begin{thmc}[{[5,8]}]
Assume that $b\in BMO(\mathbb R^n)$. Then for any $1\le p<\infty$, we have
$$\sup_B\bigg(\frac{1}{|B|}\int_B\big|b(x)-b_B\big|^p\,dx\bigg)^{1/p}\le C\|b\|_*.$$
\end{thmc}
Next we shall define the weighted Morrey space and give one of the results relevant to this paper. For further details, we refer the readers to [9].
\newtheorem{def1}{Definition}
\begin{def1}
Let $1\le p<\infty$, $0<\kappa<1$ and $w$ be a weight function. Then the weighted Morrey space is defined by
$$L^{p,\kappa}(w)=\{f\in L^p_{loc}(w):\|f\|_{L^{p,\kappa}(w)}<\infty\},$$
where
$$\|f\|_{L^{p,\kappa}(w)}=\sup_{B}\left(\frac{1}{w(B)^\kappa}\int_B|f(x)|^pw(x)\,dx\right)^{1/p}$$
and the supremum is taken over all balls $B$ in $\mathbb R^n$.
\end{def1}
In [9], the authors obtained the following result.
\newtheorem*{thmd}{Theorem D}
\begin{thmd}
If $1<p<\infty$, $0<\kappa<1$ and $w\in A_p$, then the Hardy-Littlewood maximal operator $M$ is bounded on $L^{p,\kappa}(w)$.
\end{thmd}
We are going to conclude this section by giving several results concerning the weighted boundedness of rough operators $M_\Omega$, $T_\Omega$ and $\mu_\Omega$ on the weighted $L^p$ spaces. Given a Muckenhoupt's weight function $w$ on $\mathbb R^n$, for $1\le p<\infty$, we denote by $L^p_w(\mathbb R^n)$ the space of all functions satisfying
$$\|f\|_{L^p_w}=\left(\int_{\mathbb R^n}|f(x)|^pw(x)\,dx\right)^{1/p}<\infty.$$
\newtheorem*{thme}{Theorem E}
\begin{thme}[{[4]}]
Suppose that $\Omega\in L^q(S^{n-1})$, $1<q<\infty$. Then for every $q'\le p<\infty$ and $w\in A_{p/{q'}}$, there is a constant $C$ independent of $f$ such that
$$\|M_\Omega(f)\|_{L^p_w}\le C\|f\|_{L^p_w}$$
$$\|T_\Omega(f)\|_{L^p_w}\le C\|f\|_{L^p_w}.$$
\end{thme}
\newtheorem*{thmf}{Theorem F}
\begin{thmf}[{[2]}]
Suppose that $\Omega\in L^q(S^{n-1})$, $1<q\le\infty$. Then for every $q'<p<\infty$ and $w\in A_{p/{q'}}$, there is a constant $C$ independent of $f$ such that
$$\|\mu_\Omega(f)\|_{L^p_w}\le C\|f\|_{L^p_w}.$$
\end{thmf}
\newtheorem*{thmg}{Theorem G}
\begin{thmg}[{[3]}]
Suppose that $\Omega\in L^q(S^{n-1})$ with $1<q\le\infty$ and $b\in BMO(\mathbb R^n)$. Then for $q'<p<\infty$ and $w\in A_{p/{q'}}$, there is a constant $C>0$ independent of $f$ such that
$$\big\|[b,\mu_\Omega](f)\big\|_{L^p_w}\le C\|f\|_{L^p_w}.$$
\end{thmg}
Throughout this article, we will use $C$ to denote a positive constant, which is independent of the main parameters and not necessarily the same at each occurrence. By $A\sim B$, we mean that there exists a constant $C>1$ such that $\frac1C\le\frac AB\le C$.
\newpage
\noindent\textbf{\large{3. Proof of Theorem 1}}
\par
First, by using H\"older's inequality, we can easily see that
$$M_\Omega f(x)\le C\cdot\|\Omega\|_{L^q(S^{n-1})}M_{q'}(f)(x),$$
where $M_{q'}(f)(x)=M(|f|^{q'})(x)^{1/{q'}}.$ Then for $q'<p<\infty$ and $w\in A_{p/{q'}}$, it follows immediately from Theorem D that
\begin{equation*}
\begin{split}
\big\|M_{q'}(f)\big\|_{L^{p,\kappa}(w)}&=\big\|M(|f|^{q'})\big\|_{L^{p/{q'},\kappa}(w)}^{1/{q'}}\\
&\le C\big\||f|^{q'}\big\|_{L^{p/{q'},\kappa}(w)}^{1/{q'}}\\
&\le C\|f\|_{L^{p,\kappa}(w)}.
\end{split}
\end{equation*}
\par
Now we consider the case $p=q'$. Fix a ball $B=B(x_0,r_B)\subseteq\mathbb R^n$ and decompose $f=f_1+f_2$, where $f_1=f\chi_{_{2B}}$, $\chi_{_{2B}}$ denotes the characteristic function of $2B$. Since $M_\Omega$ is a sublinear operator, then we have
\begin{equation*}
\begin{split}
&\frac{1}{w(B)^{\kappa/p}}\Big(\int_B|M_\Omega f(x)|^pw(x)\,dx\Big)^{1/p}\\
\le\,&\frac{1}{w(B)^{\kappa/p}}\Big(\int_B|M_\Omega f_1(x)|^pw(x)\,dx\Big)^{1/p}\\
&+\frac{1}{w(B)^{\kappa/p}}\Big(\int_B|M_\Omega f_2(x)|^pw(x)\,dx\Big)^{1/p}\\
=\,&I_1+I_2.
\end{split}
\end{equation*}
Theorem E and Lemma A imply
\begin{equation}
\begin{split}
I_1&\le C\cdot\frac{1}{w(B)^{\kappa/p}}\Big(\int_{2B}|f(x)|^pw(x)\,dx\Big)^{1/p}\\
&\le C\|f\|_{L^{p,\kappa}(w)}\cdot\frac{w(2B)^{\kappa/p}}{w(B)^{\kappa/p}}\\
&\le C\|f\|_{L^{p,\kappa}(w)}.
\end{split}
\end{equation}
We turn to estimate the term $I_2$. For any given $r>0$ and $x\in B$, by H\"older's inequality and the $A_1$ condition, we thus obtain
\begin{equation*}
\begin{split}
&\frac{1}{r^n}\int_{|y|<r}|\Omega(y')f_2(x-y)|\,dy\\
\le\,&\frac{1}{r^n}\Big(\int_{|y|<r}\big|\Omega(y')\big|^q\,dy\Big)^{1/q}\Big(\int_{|y|<r}\big|f_2(x-y)\big|^p\,dy\Big)^{1/p}\\
\le\,&C\cdot\|\Omega\|_{L^q(S^{n-1})}\bigg(\frac{1}{|B(x,r)|}\int_{B(x,r)}\big|f_2(y)\big|^p\,dy\bigg)^{1/p}\\
\end{split}
\end{equation*}
\begin{equation*}
\begin{split}
\le\,&C\cdot\|\Omega\|_{L^q(S^{n-1})}\bigg(\frac{1}{w(B(x,r))}\int_{B(x,r)}\big|f_2(y)\big|^pw(y)\,dy\bigg)^{1/p}.
\end{split}
\end{equation*}
A simple geometric observation shows that when $x\in B(x_0,r_B)$ and $y\in B(x,r)\cap(2B(x_0,r_B))^c$, then we have $B(x_0,r_B)\subseteq 3B(x,r).$ Hence
\begin{equation*}
\begin{split}
\frac{1}{r^n}\int_{|y|<r}|\Omega(y')f_2(x-y)|\,dy&\le C\|f\|_{L^{p,\kappa}(w)}\cdot\frac{1}{w(B(x,r))^{(1-\kappa)/p}}\\
&\le C\|f\|_{L^{p,\kappa}(w)}\cdot\frac{1}{w(B(x_0,r_B))^{(1-\kappa)/p}}.
\end{split}
\end{equation*}
Taking the supremum over all $r>0$, we can get
\begin{equation*}
\big|M_\Omega(f_2)(x)\big|\le C\|f\|_{L^{p,\kappa}(w)}\cdot\frac{1}{w(B(x_0,r_B))^{(1-\kappa)/p}},
\end{equation*}
which implies
\begin{equation}
I_2\le C\|f\|_{L^{p,\kappa}(w)}.
\end{equation}
Combining the above inequality (2) with (1) and taking the supremum over all balls $B\subseteq\mathbb R^n$, we obtain the desired result.\\
\textbf{\large{4. Proofs of Theorems 2 and 3}}
\begin{proof}[Proof of Theorem 2]
Fix a ball $B=B(x_0,r_B)$ and decompose $f=f_1+f_2$, where $f_1=f\chi_{_{2B}}$. Then we have
\begin{equation*}
\begin{split}
&\frac{1}{w(B)^{\kappa/p}}\Big(\int_B|T_\Omega f(x)|^pw(x)\,dx\Big)^{1/p}\\
\le\,&\frac{1}{w(B)^{\kappa/p}}\Big(\int_B|T_\Omega f_1(x)|^pw(x)\,dx\Big)^{1/p}\\
&+\frac{1}{w(B)^{\kappa/p}}\Big(\int_B|T_\Omega f_2(x)|^pw(x)\,dx\Big)^{1/p}\\
=\,&J_1+J_2.
\end{split}
\end{equation*}
Theorem E and Lemma A give
\begin{equation*}
\begin{split}
J_1&\le C\cdot\frac{1}{w(B)^{\kappa/p}}\Big(\int_{2B}|f(x)|^pw(x)\,dx\Big)^{1/p}\\
&\le C\|f\|_{L^{p,\kappa}(w)}\cdot\frac{w(2B)^{\kappa/p}}{w(B)^{\kappa/p}}\\
&\le C\|f\|_{L^{p,\kappa}(w)}.
\end{split}
\end{equation*}
In order to estimate $J_2$, we first deduce from H\"older's inequality that
\begin{equation*}
\begin{split}
&|T_\Omega(f_2)(x)|\\
=\,&\bigg|\int_{(2B)^c}\frac{\Omega((x-y)')}{|x-y|^n}f(y)\,dy\bigg|\\
\le\,&\sum_{j=1}^\infty\Big(\int_{2^{j+1}B\backslash2^jB}\big|\Omega((x-y)')\big|^q\,dy\Big)^{1/q}\Big(\int_{2^{j+1}B\backslash2^jB}\frac{|f(y)|^{q'}}{|x-y|^{nq'}}\,dy\Big)^{1/{q'}}.
\end{split}
\end{equation*}
When $x\in B$ and $y\in2^{j+1}B\backslash2^jB$, then by a direct calculation, we can see that $2^{j-1}r_B\le|y-x|<2^{j+2}r_B$. Hence
\begin{equation}
\Big(\int_{2^{j+1}B\backslash2^jB}\big|\Omega((x-y)')\big|^q\,dy\Big)^{1/q}\le C\cdot\|\Omega\|_{L^q(S^{n-1})}|2^{j+1}B|^{1/q}.
\end{equation}
We also note that if $x\in B$, $y\in(2B)^c$, then $|y-x|\sim|y-x_0|$. Consequently
$$\Big(\int_{2^{j+1}B\backslash2^jB}\frac{|f(y)|^{q'}}{|x-y|^{nq'}}\,dy\Big)^{1/{q'}}\le\frac{1}{|2^{j+1}B|}\Big(\int_{2^{j+1}B}|f(y)|^{q'}\,dy\Big)^{1/{q'}}.$$
So we have
$$\big|T_\Omega(f_2)(x)\big|\le C\sum_{j=1}^\infty\Big(\frac{1}{|2^{j+1}B|}\int_{2^{j+1}B}|f(y)|^{q'}\,dy\Big)^{1/{q'}}.$$
We shall consider two cases. When $p=q'$, then by the $A_1$ condition, we get
\begin{equation}
\begin{split}
\big|T_\Omega(f_2)(x)\big|&\le C\sum_{j=1}^\infty\Big(\frac{1}{w(2^{j+1}B)}\int_{2^{j+1}B}|f(y)|^{p}w(y)\,dy\Big)^{1/{p}}\\
&\le C\|f\|_{L^{p,\kappa}(w)}\sum_{j=1}^\infty\frac{1}{w(2^{j+1}B)^{(1-\kappa)/p}}.
\end{split}
\end{equation}
When $p>q'$, set $s=p/{q'}>1$. Then it follows from the H\"older's inequality and the $A_{s}$ condition that
\begin{equation}
\begin{split}
\big|T_\Omega(f_2)(x)\big|\le\,& C\sum_{j=1}^\infty\frac{1}{|2^{j+1}B|^{1/{q'}}}\Big(\int_{2^{j+1}B}|f(y)|^{p}w(y)\,dy\Big)^{1/{p}}\\
&\times\Big(\int_{2^{j+1}B}w^{-s'/s}(y)\,dy\Big)^{1/{s'q'}}\\
\le\,&C\sum_{j=1}^\infty\Big(\frac{1}{w(2^{j+1}B)}\int_{2^{j+1}B}|f(y)|^{p}w(y)\,dy\Big)^{1/{p}}\\
\le\,& C\|f\|_{L^{p,\kappa}(w)}\sum_{j=1}^\infty\frac{1}{w(2^{j+1}B)^{(1-\kappa)/p}}.
\end{split}
\end{equation}
Hence, for every $q'\le p<\infty$, by the estimates (4) and (5), we obtain
$$J_2\le C\|f\|_{L^{p,\kappa}(w)}\sum_{j=1}^\infty\left(\frac{w(B)}{w(2^{j+1}B)}\right)^{(1-\kappa)/p}.$$
Since $w\in A_{p/{q'}}$, then there exists $r>1$ such that $w\in RH_r$. By using Lemma B, we thus get
\begin{equation}
\frac{w(B)}{w(2^{j+1}B)}\le C\left(\frac{|B|}{|2^{j+1}B|}\right)^{(r-1)/r}.
\end{equation}
Therefore
\begin{equation*}
\begin{split}
J_2&\le C\|f\|_{L^{p,\kappa}(w)}\sum_{j=1}^\infty\Big(\frac{1}{2^{jn}}\Big)^{{(1-\kappa)(r-1)}/{pr}}\\
&\le C\|f\|_{L^{p,\kappa}(w)},
\end{split}
\end{equation*}
where the last series is convergent since ${(1-\kappa)(r-1)}/{pr}>0$. Using the estimates for $J_1$ and $J_2$ and taking the supremum over all balls $B\subseteq\mathbb R^n$, we complete the proof of Theorem 2.
\end{proof}
\begin{proof}[Proof of Theorem 3]
As in the proof of Theorem 2, we can write
\begin{equation*}
\begin{split}
&\frac{1}{w(B)^{\kappa/p}}\Big(\int_B\big|[b,T_\Omega]f(x)\big|^pw(x)\,dx\Big)^{1/p}\\
\le\,&\frac{1}{w(B)^{\kappa/p}}\Big(\int_B\big|[b,T_\Omega]f_1(x)\big|^pw(x)\,dx\Big)^{1/p}\\
&+\frac{1}{w(B)^{\kappa/p}}\Big(\int_B\big|[b,T_\Omega]f_2(x)\big|^pw(x)\,dx\Big)^{1/p}\\
=\,&J'_1+J'_2.
\end{split}
\end{equation*}
By Theorem E and the well-known boundedness criterion for the commutators of linear operators, which was obtained by Alvarez, Bagby, Kurtz and P\'erez(see [1]), we see that $[b,T_\Omega]$ is bounded on $L^p_w$ for all $q'<p<\infty$ and $w\in A_{p/{q'}}$. This together with Lemma A yield
\begin{equation}
\begin{split}
J'_1&\le C\|b\|_*\cdot\frac{1}{w(B)^{\kappa/p}}\Big(\int_{2B}|f(x)|^pw(x)\,dx\Big)^{1/p}\\
&\le C\|b\|_*\|f\|_{L^{p,\kappa}(w)}\cdot\frac{w(2B)^{\kappa/p}}{w(B)^{\kappa/p}}\\
&\le C\|b\|_*\|f\|_{L^{p,\kappa}(w)}.
\end{split}
\end{equation}
We now turn to deal with the term $J'_2$. For any given $x\in B$, we have
\begin{equation*}
\begin{split}
\big|[b,T_\Omega]f_2(x)\big|\le&\, \big|b(x)-b_B\big|\cdot\int_{(2B)^c}\frac{|\Omega((x-y)')|}{|x-y|^n}|f(y)|\,dy\\
&+\int_{(2B)^c}\frac{|\Omega((x-y)')|}{|x-y|^n}|b(y)-b_B||f(y)|\,dy\\
=\,&\mbox{\upshape I+II}.
\end{split}
\end{equation*}
In the proof of Theorem 2, for any $q'<p<\infty$, we have already showed
\begin{equation*}
\mbox{\upshape I}\le C|b(x)-b_B|\cdot\|f\|_{L^{p,\kappa}(w)}\sum_{j=1}^\infty\frac{1}{w(2^{j+1}B)^{(1-\kappa)/p}}.
\end{equation*}
Consequently
\begin{equation*}
\begin{split}
&\frac{1}{w(B)^{\kappa/p}}\Big(\int_B \mbox{\upshape I}^p\,w(x)\,dx\Big)^{1/p}\\
\le \,&C\|f\|_{L^{p,\kappa}(w)}\frac{1}{w(B)^{\kappa/p}}\cdot\sum_{j=1}^\infty \frac{1}{w(2^{j+1}B)^{(1-\kappa)/p}} \cdot\Big(\int_B\big|b(x)-b_B\big|^pw(x)\,dx\Big)^{1/p}\\
=\,&C\|f\|_{L^{p,\kappa}(w)}\sum_{j=1}^\infty\frac{w(B)^{(1-\kappa)/p}}{w(2^{j+1}B)^{(1-\kappa)/p}}
\cdot\Big(\frac{1}{w(B)}\int_B\big|b(x)-b_B\big|^pw(x)\,dx\Big)^{1/p}.
\end{split}
\end{equation*}
Using the same arguments as that of Theorem 2, we can see that the above summation is bounded by a constant. Hence
\begin{equation*}
\begin{split}
\frac{1}{w(B)^{\kappa/p}}\Big(\int_B \mbox{\upshape I}^p\,w(x)\,dx\Big)^{1/p}\le & C\|f\|_{L^{p,\kappa}(w)}\Big(\frac{1}{w(B)}\int_B\big|b(x)-b_B\big|^pw(x)\,dx\Big)^{1/p}.
\end{split}
\end{equation*}
Since $w\in A_{p/{q'}}$, then $w\in A_p$. As before, there exists a number $r>1$ such that $w\in RH_r$. By the reverse H\"older's inequality and Theorem C, we get
\begin{equation}
\begin{split}
&\Big(\frac{1}{w(B)}\int_B\big|b(x)-b_B\big|^pw(x)\,dx\Big)^{1/p}\\
\le\, &C\cdot\frac{1}{w(B)^{1/p}}\Big(\int_B\big|b(x)-b_B\big|^{pr'}\,dx\Big)^{1/{pr'}}\Big(\int_Bw(x)^r\,dx\Big)^{1/{pr}}\\
\le\, &C\cdot\Big(\frac{1}{|B|}\int_B\big|b(x)-b_B\big|^{pr'}\,dx\Big)^{1/{pr'}}\\
\le\, &C\|b\|_*.
\end{split}
\end{equation}
So we have
\begin{equation}
\frac{1}{w(B)^{\kappa/p}}\Big(\int_B \mbox{\upshape I}^p\,w(x)\,dx\Big)^{1/p}\le C\|b\|_*\|f\|_{L^{p,\kappa}(w)}.
\end{equation}
On the other hand, it follows from H\"older's inequality and (3) that
\begin{equation*}
\begin{split}
\mbox{\upshape II}\le\, & C\sum_{j=1}^\infty\Big(\frac{1}{|2^{j+1}B|}\int_{2^{j+1}B}\big|b(y)-b_B\big|^{q'}\big|f(y)\big|^{q'}\,dy\Big)^{1/{q'}}\\
\le\,&C\sum_{j=1}^\infty\Big(\frac{1}{|2^{j+1}B|}\int_{2^{j+1}B}\big|b(y)-b_{2^{j+1}B}\big|^{q'}\big|f(y)\big|^{q'}\,dy\Big)^{1/{q'}}\\
&+C\sum_{j=1}^\infty\big|b_{2^{j+1}B}-b_B\big|\cdot\Big(\frac{1}{|2^{j+1}B|}\int_{2^{j+1}B}\big|f(y)\big|^{q'}\,dy\Big)^{1/{q'}}\\
=\,&\mbox{\upshape III+IV}.
\end{split}
\end{equation*}
Set $s=p/{q'}>1$. Then by using H\"older's inequality, we thus obtain
\begin{align}
&\Big(\int_{2^{j+1}B}\big|b(y)-b_{2^{j+1}B}\big|^{q'}\big|f(y)\big|^{q'}\,dy\Big)^{1/{q'}}\\
\le\,&\Big(\int_{2^{j+1}B}\big|b(y)-b_{2^{j+1}B}\big|^{q's'}w^{-{s'}/s}(y)\,dy\Big)^{1/{q's'}}
\Big(\int_{2^{j+1}B}\big|f(y)\big|^pw(y)\,dy\Big)^{1/p}\notag\\
\le\,&C\|f\|_{L^{p,\kappa}(w)}\cdot w(2^{j+1}B)^{\kappa/p}
\Big(\int_{2^{j+1}B}\big|b(y)-b_{2^{j+1}B}\big|^{q's'}w^{-{s'}/s}(y)\,dy\Big)^{1/{q's'}}\notag.
\end{align}
Let $v(y)=w^{-{s'}/s}(y)=w^{1-s'}(y)$. Then we have $v\in A_{s'}$ because $w\in A_s$(see [6]), which implies $v\in A_{q's'}$. Following along the same lines as that of (8), we can get
\begin{equation}
\Big(\frac{1}{v(2^{j+1}B)}\int_{2^{j+1}B}\big|b(y)-b_{2^{j+1}B}\big|^{q's'}v(y)\,dy\Big)^{1/{q's'}}\le C\|b\|_*.
\end{equation}
Substituting the above inequality (11) into (10), we thus have
\begin{equation}
\begin{split}
&\Big(\int_{2^{j+1}B}\big|b(y)-b_{2^{j+1}B}\big|^{q'}\big|f(y)\big|^{q'}\,dy\Big)^{1/{q'}}\\
\le& C\|b\|_*\|f\|_{L^{p,\kappa}(w)}\cdot w(2^{j+1}B)^{\kappa/p}\cdot v(2^{j+1}B)^{1/{q's'}}\\
\le& C\|b\|_*\|f\|_{L^{p,\kappa}(w)}\cdot|2^{j+1}B|^{1/{q'}}\cdot w(2^{j+1}B)^{(\kappa-1)/p}.
\end{split}
\end{equation}
Hence
\begin{align}
\frac{1}{w(B)^{\kappa/p}}\Big(\int_B \mbox{\upshape III}^p\,w(x)\,dx\Big)^{1/p}&\le C\|b\|_*\|f\|_{L^{p,\kappa}(w)}\sum_{j=1}^\infty\frac{w(B)^{(1-\kappa)/p}}{w(2^{j+1}B)^{(1-\kappa)/p}}\notag\\
&\le C\|b\|_*\|f\|_{L^{p,\kappa}(w)}.
\end{align}
Now let's deal with the last term \mbox{\upshape IV}. Since $b\in BMO(\mathbb R^n)$, then a simple computation shows that
\begin{equation}
|b_{2^{j+1}B}-b_B|\le C\cdot j\|b\|_*.
\end{equation}
It follows immediately from the inequalities (5) and (14) that
\begin{equation*}
\begin{split}
\mbox{\upshape IV}&\le C\|b\|_*\sum_{j=1}^\infty j\cdot\Big(\frac{1}{|2^{j+1}B|}\int_{2^{j+1}B}\big|f(y)\big|^{q'}\,dy\Big)^{1/{q'}}\\
&\le C\|b\|_*\|f\|_{L^{p,\kappa}(w)}\sum_{j=1}^\infty j\cdot w(2^{j+1}B)^{(\kappa-1)/p}.
\end{split}
\end{equation*}
Therefore, by the estimate (6), we obtain
\begin{align}
\frac{1}{w(B)^{\kappa/p}}\Big(\int_B \mbox{\upshape IV}^p\,w(x)\,dx\Big)^{1/p}&\le C\|b\|_*\|f\|_{L^{p,\kappa}(w)}\sum_{j=1}^\infty j\cdot\frac{w(B)^{(1-\kappa)/p}}{w(2^{j+1}B)^{(1-\kappa)/p}}\notag\\
&\le C\|b\|_*\|f\|_{L^{p,\kappa}(w)}\sum_{j=1}^\infty\frac{j}{2^{jn\theta}}\notag\\
&\le C\|b\|_*\|f\|_{L^{p,\kappa}(w)},
\end{align}
where $w\in RH_r$ and $\theta=(1-\kappa)(r-1)/{pr}$. Summarizing the estimates (13) and (15) derived above, we can get
\begin{equation}
\frac{1}{w(B)^{\kappa/p}}\Big(\int_B \mbox{\upshape II}^p\,w(x)\,dx\Big)^{1/p}\le C\|b\|_*\|f\|_{L^{p,\kappa}(w)}.
\end{equation}
Combining the inequalities (7), (9) with the inequality (16) and taking the supremum over all balls $B\subseteq\mathbb R^n$, we conclude the proof of Theorem 3.
\end{proof}
\noindent\textbf{\large{5. Proofs of Theorems 4 and 5}}
\begin{proof}[Proof of Theorem 4]
Fix a ball $B=B(x_0,r_B)\subseteq\mathbb R^n$. Let $f=f_1+f_2$, where $f_1=f\chi_{_{2B}}$. Then we have
\begin{equation*}
\begin{split}
&\frac{1}{w(B)^{\kappa/p}}\Big(\int_B|\mu_\Omega f(x)|^pw(x)\,dx\Big)^{1/p}\\
\le\,&\frac{1}{w(B)^{\kappa/p}}\Big(\int_B|\mu_\Omega f_1(x)|^pw(x)\,dx\Big)^{1/p}\\
&+\frac{1}{w(B)^{\kappa/p}}\Big(\int_B|\mu_\Omega f_2(x)|^pw(x)\,dx\Big)^{1/p}\\
=\,&K_1+K_2.
\end{split}
\end{equation*}
Theorem F and Lemma A imply
\begin{equation*}
\begin{split}
K_1&\le C\cdot\frac{1}{w(B)^{\kappa/p}}\Big(\int_{2B}|f(x)|^pw(x)\,dx\Big)^{1/p}\\
&\le C\|f\|_{L^{p,\kappa}(w)}\cdot\frac{w(2B)^{\kappa/p}}{w(B)^{\kappa/p}}\\
&\le C\|f\|_{L^{p,\kappa}(w)}.
\end{split}
\end{equation*}
To estimate $K_2$, observe that when $x\in B$ and $y\in 2^{j+1}B\backslash2^j B$($j\ge1$), then
$$t\ge|x-y|\ge|y-x_0|-|x-x_0|\ge2^{j-1}r_B.$$
Therefore
\begin{equation*}
\begin{split}
\big|\mu_\Omega(f_2)(x)\big|&=\left(\int_0^\infty\Big|\int_{(2B)^c\cap\{y:|x-y|\le t\}}\frac{\Omega(x-y)}{|x-y|^{n-1}}f(y)\,dy\Big|^2\frac{dt}{t^3}\right)^{1/2}\\
&\le\sum_{j=1}^\infty\Big(\int_{2^{j+1}B\backslash2^j B}\frac{|\Omega(x-y)|}{|x-y|^{n-1}}|f(y)|\,dy\Big)\cdot\Big(\int_{2^{j-1}r_B}^\infty\frac{dt}{t^3}\Big)^{1/2}\\
&\le C\sum_{j=1}^\infty\frac{1}{|2^{j+1}B|^{1/n}}\cdot\int_{2^{j+1}B\backslash2^j B}\frac{|\Omega(x-y)|}{|x-y|^{n-1}}|f(y)|\,dy.
\end{split}
\end{equation*}
When $\Omega\in L^\infty(S^{n-1})$, then by assumption, we have $w\in A_p$, $1<p<\infty$. It follows from the H\"older's inequality and the $A_p$ condition that
\begin{align}
\big|\mu_\Omega(f_2)(x)\big|\le\,& C\|\Omega\|_{L^\infty(S^{n-1})}\sum_{j=1}^\infty\frac{1}{|2^{j+1}B|^{1/n}}\cdot\frac{1}{|2^{j+1}B|^{(n-1)/n}}\int_{2^{j+1}B}|f(y)|\,dy\notag\\
\le\, & C\|\Omega\|_{L^\infty(S^{n-1})}\sum_{j=1}^\infty\frac{1}{|2^{j+1}B|}\Big(\int_{2^{j+1}B}|f(y)|^pw(y)\,dy\Big)^{1/p}\notag\\
&\times\Big(\int_{2^{j+1}B}w(y)^{-p'/p}\,dy\Big)^{1/{p'}}\\
\le\, &C\|\Omega\|_{L^\infty(S^{n-1})}\|f\|_{L^{p,\kappa}(w)}\sum_{j=1}^\infty w(2^{j+1}B)^{(\kappa-1)/p}\notag.
\end{align}
When $\Omega\in L^q(S^{n-1})$, $1<q<\infty$, by the inequalities (3) and (5), we get
\begin{equation}
\begin{split}
\big|\mu_\Omega(f_2)(x)\big|&\le C\|\Omega\|_{L^q(S^{n-1})}\sum_{j=1}^\infty\Big(\frac{1}{|2^{j+1}B|}\int_{2^{j+1}B}|f(y)|^{q'}\,dy\Big)^{1/{q'}}\\
&\le C\|\Omega\|_{L^q(S^{n-1})}\|f\|_{L^{p,\kappa}(w)}\sum_{j=1}^\infty w(2^{j+1}B)^{(\kappa-1)/p}.
\end{split}
\end{equation}
Hence, for $1<q\le\infty$, $q'<p<\infty$, by the estimates (17) and (18), we have
\begin{equation*}
\begin{split}
K_2&\le C\|f\|_{L^{p,\kappa}(w)}\sum_{j=1}^\infty\left(\frac{w(B)}{w(2^{j+1}B)}\right)^{(1-\kappa)/p}\\
&\le C\|f\|_{L^{p,\kappa}(w)}.
\end{split}
\end{equation*}
Using the above estimates for $K_1$ and $K_2$ and taking the supremum over all balls $B\subseteq\mathbb R^n$, we get our desired result.
\end{proof}
\begin{proof}[Proof of Theorem 5]
As before, we can write
\begin{equation*}
\begin{split}
&\frac{1}{w(B)^{\kappa/p}}\Big(\int_B\big|[b,\mu_\Omega]f(x)\big|^pw(x)\,dx\Big)^{1/p}\\
\le\,&\frac{1}{w(B)^{\kappa/p}}\Big(\int_B\big|[b,\mu_\Omega]f_1(x)\big|^pw(x)\,dx\Big)^{1/p}\\
&+\frac{1}{w(B)^{\kappa/p}}\Big(\int_B\big|[b,\mu_\Omega]f_2(x)\big|^pw(x)\,dx\Big)^{1/p}\\
=\,&K'_1+K'_2.
\end{split}
\end{equation*}
Theorem G and Lemma A yield
\begin{equation*}
\begin{split}
K'_1&\le C\|b\|_*\|f\|_{L^{p,\kappa}(w)}\frac{w(2B)^{\kappa/p}}{w(B)^{\kappa/p}}\\
&\le C\|b\|_*\|f\|_{L^{p,\kappa}(w)}.
\end{split}
\end{equation*}
Finally, let us deal with the term $K'_2$. For any fixed $x\in B$, we have
\begin{equation*}
\begin{split}
&\big|[b,\mu_\Omega]f_2(x)\big|\\
\le\,&|b(x)-b_B|\left(\int_0^\infty\Big|\int_{(2B)^c\cap\{y:|x-y|\le t\}}\frac{\Omega(x-y)}{|x-y|^{n-1}}f(y)\,dy\Big|^2\frac{dt}{t^3}\right)^{1/2}\\
&+\left(\int_0^\infty\Big|\int_{(2B)^c\cap\{y:|x-y|\le t\}}\frac{\Omega(x-y)}{|x-y|^{n-1}}[b(y)-b_B]f(y)\,dy\Big|^2\frac{dt}{t^3}\right)^{1/2}\\
=\,&\mbox{\upshape I+II}.
\end{split}
\end{equation*}
In the proof of Theorem 4, for any $q'<p<\infty$, we have already proved
\begin{equation*}
\mbox{\upshape I}\le C|b(x)-b_B|\cdot\|f\|_{L^{p,\kappa}(w)}\sum_{j=1}^\infty w(2^{j+1}B)^{(\kappa-1)/p}.
\end{equation*}
Following the same lines as in the proof of Theorem 3, we obtain
\begin{equation*}
\frac{1}{w(B)^{\kappa/p}}\Big(\int_B \mbox{\upshape I}^p\,w(x)\,dx\Big)^{1/p}\le C\|b\|_*\|f\|_{L^{p,\kappa}(w)}.
\end{equation*}
On the other hand, we note that when $x\in B$ and $y\in 2^{j+1}B\backslash2^j B$($j\ge1$), then we have $t\ge2^{j-1}r_B$. Consequently
\begin{equation*}
\begin{split}
\mbox{\upshape II}\le\, &C\sum_{j=1}^\infty\frac{1}{|2^{j+1}B|^{1/n}}\cdot\int_{2^{j+1}B\backslash2^j B}\frac{|\Omega(x-y)|}{|x-y|^{n-1}}|b(y)-b_B||f(y)|\,dy\\
\le\,&C\sum_{j=1}^\infty\frac{1}{|2^{j+1}B|^{1/n}}\cdot\int_{2^{j+1}B\backslash2^j B}\frac{|\Omega(x-y)|}{|x-y|^{n-1}}|b(y)-b_{2^{j+1}B}||f(y)|\,dy\\
&+C\sum_{j=1}^\infty\frac{|b_{2^{j+1}B}-b_B|}{|2^{j+1}B|^{1/n}}\cdot\int_{2^{j+1}B\backslash2^j B}\frac{|\Omega(x-y)|}{|x-y|^{n-1}}|f(y)|\,dy\\
=\,&\mbox{\upshape III+IV}.
\end{split}
\end{equation*}
When $\Omega\in L^\infty(S^{n-1})$, then it follows from H\"older's inequality that
\begin{equation*}
\begin{split}
\mbox{\upshape III}\le \,& C\|\Omega\|_{L^\infty(S^{n-1})}\sum_{j=1}^\infty\frac{1}{|2^{j+1}B|}\int_{2^{j+1}B}|b(y)-b_{2^{j+1}B}||f(y)|\,dy\\
\le\,&C\|\Omega\|_{L^\infty(S^{n-1})}\sum_{j=1}^\infty\frac{1}{|2^{j+1}B|}\Big(\int_{2^{j+1}B}\big|b(y)-b_{2^{j+1}B}\big|^{p'}w^{-{p'}/p}(y)\,dy\Big)^{1/{p'}}\\
&\times\Big(\int_{2^{j+1}B}\big|f(y)\big|^pw(y)\,dy\Big)^{1/p}\\
\le\,&C\|\Omega\|_{L^\infty(S^{n-1})}\|f\|_{L^{p,\kappa}(w)}\sum_{j=1}^\infty \frac{1}{|2^{j+1}B|}\cdot w(2^{j+1}B)^{\kappa/p}\\
&\times\Big(\int_{2^{j+1}B}\big|b(y)-b_{2^{j+1}B}\big|^{p'}w^{-{p'}/p}(y)\,dy\Big)^{1/{p'}}.
\end{split}
\end{equation*}
Set $u(y)=w^{-{p'}/p}(y)=w^{1-p'}(y)$. In this case, since $w\in A_p$, then we have $u\in A_{p'}$, it follows from the inequality (8) and the $A_p$ condition that
\begin{align}
\mbox{\upshape III}&\le C\|\Omega\|_{L^\infty(S^{n-1})}\|b\|_*\|f\|_{L^{p,\kappa}(w)}\cdot\sum_{j=1}^\infty \frac{1}{|2^{j+1}B|} w(2^{j+1}B)^{\kappa/p}\cdot u(2^{j+1}B)^{1/{p'}}\notag\\
&\le C\|\Omega\|_{L^\infty(S^{n-1})}\|b\|_*\|f\|_{L^{p,\kappa}(w)}\cdot\sum_{j=1}^\infty w(2^{j+1}B)^{(\kappa-1)/p}.
\end{align}
When $\Omega\in L^q(S^{n-1})$, then by using H\"older's inequality, the inequalities (3) and (12), we can deduce
\begin{align}
\mbox{\upshape III}&\le C\|\Omega\|_{L^q(S^{n-1})}\sum_{j=1}^\infty\frac{1}{|2^{j+1}B|^{1/{q'}}}\Big(\int_{2^{j+1}B}\big|b(y)-b_{2^{j+1}B}\big|^{q'}\big|f(y)\big|^{q'}\,dy\Big)^{1/{q'}}\notag\\
&\le C\|\Omega\|_{L^q(S^{n-1})}\|b\|_*\|f\|_{L^{p,\kappa}(w)}\cdot\sum_{j=1}^\infty w(2^{j+1}B)^{(\kappa-1)/p}.
\end{align}
Hence, for $1<q\le\infty$, $q'<p<\infty$, by the estimates (19) and (20), we get
\begin{align}
\frac{1}{w(B)^{\kappa/p}}\Big(\int_B \mbox{\upshape III}^p\,w(x)\,dx\Big)^{1/p}&\le C\|b\|_*\|f\|_{L^{p,\kappa}(w)}\sum_{j=1}^\infty\frac{w(B)^{(1-\kappa)/p}}{w(2^{j+1}B)^{(1-\kappa)/p}}\notag\\
&\le C\|b\|_*\|f\|_{L^{p,\kappa}(w)}\notag.
\end{align}
Again, in Theorem 4, we have already obtained the following inequality
\begin{equation*}
\frac{1}{|2^{j+1}B|^{1/n}}\cdot\int_{2^{j+1}B\backslash2^j B}\frac{|\Omega(x-y)|}{|x-y|^{n-1}}|f(y)|\,dy\le C\|f\|_{L^{p,\kappa}(w)}\cdot w(2^{j+1}B)^{(\kappa-1)/p}.
\end{equation*}
From (14), it follows immediately that
\begin{equation*}
\mbox{\upshape IV}\le C\|b\|_*\|f\|_{L^{p,\kappa}(w)}\sum_{j=1}^\infty j\cdot w(2^{j+1}B)^{(\kappa-1)/p}.
\end{equation*}
The rest of the proof is exactly the same as that of (15), we finally obtain
\begin{equation*}
\frac{1}{w(B)^{\kappa/p}}\Big(\int_B \mbox{\upshape IV}^p\,w(x)\,dx\Big)^{1/p}\le C\|b\|_*\|f\|_{L^{p,\kappa}(w)}.
\end{equation*}
Therefore, by combining the above estimates and taking the supremum over all balls $B\subseteq\mathbb R^n$, we conclude the proof of Theorem 5.
\end{proof}

\end{document}